\documentclass{article}
\usepackage[letterpaper,top=2cm,bottom=2cm,left=2.5cm,right=2.5cm,marginparwidth=1.75cm]{geometry}

\usepackage{enumerate}
\usepackage{mathtools}
\usepackage{amsmath}
\usepackage{amssymb}
\usepackage{graphicx}
\usepackage{amsthm}
\usepackage[square,numbers]{natbib}
\makeatletter
\renewcommand{\@seccntformat}[1]{\csname the#1\endcsname\ }
\makeatother

\newtheorem{theorem}{Theorem}
\newtheorem{lemma}{Lemma}
\newtheorem{corollary}{Corollary}[theorem]

\theoremstyle{definition}

\newtheorem{remark}{Remark}

\title{Raabe's Formula For Gamma Function Via Riemann-Liouville Fractional Integrals And Generalized Glaisher Constants}
\author{Efe Gürel$^{1}$\\
        \small $^{1}$TÜBİTAK Natural Sciences High School, Kocaeli 41400, Turkey}
\date{}

\begin{document}

\maketitle
\begin{abstract}    
    In this paper, we prove Raabe-type integral formulas for gamma function via left and right sided Riemann-Liouville fractional integrals. As corollaries, we give the left and right sided repeated integration formulas for the log-gamma and related functions. The relationship between the generalized Glaisher constants and aforementioned integrals are investigated.\\

     
\end{abstract}

\section{Introduction}
Euler's gamma function \cite{HigherTrans} is defined for $\mathfrak{Re}(z)>0$ as the following integral
\begin{align*}
    \Gamma(z)=\int_{0}^{\infty}t^{z-1}e^{-t}dt.
\end{align*}
Gamma function satisfies the celebrated functional equation $\Gamma(z+1)=z\Gamma(z)$ and the Euler reflection formula
\begin{align*}
    \Gamma(z)\Gamma(1-z)=\frac{\pi}{\sin\pi z}, \qquad \left(z\in\mathbb{C}\setminus \mathbb{Z}_{\le 0} \right).
\end{align*}
Here and throughout, $\mathbb{C}$ and $\mathbb{Z}_{\le 0}$ denote the sets of complex numbers and non-positive integers, respectively. The digamma function, $\psi(s)$, is defined as the logarithmic derivative of the gamma function and can be expressed through the series representation for $z\in\mathbb{C}\setminus \mathbb{Z}_{\le 0}$,
\begin{align*}
    \psi(z)=\frac{d }{dz}\log\Gamma(z)=-\gamma+\sum_{n=0}^{\infty}\frac{1}{n+1}-\frac{1}{n+z}
\end{align*}
where $\gamma$ is the Euler-Mascheroni constant given by
\begin{align*}
    \gamma=\lim_{n \to \infty}\left( \sum_{k=1}^{n}\frac{1}{k}-\log n \right)=0.57721 56649 01532\ldots.
\end{align*}
The Bernoulli polynomials $B_n(x)$ \cite{StegunHandbook} are represented by the generating function  
\begin{align*}
    \frac{te^{xt}}{e^t-1}=\sum_{n=0}^{\infty}B_n(x)\frac{t^n}{n!}, \qquad(\left| t \right|<2\pi),
\end{align*}
and Bernoulli numbers are defined as $B_n=B_n(0)$. Special values of the Bernoulli numbers are $B_0=1$, $B_1=-\frac{1}{2}=-B_1(1)$, $B_n=0$ for every odd $n$ and $B_n(0)=B_n(1)$ for every $n>1$. The Riemann zeta function is defined for $\mathfrak{Re}(s)>1$ by the series (see e.g.,\cite{TheoryofZeta})
\begin{align*}
    \zeta(s)=\sum_{n=1}^{\infty}\frac{1}{n^s}
\end{align*}
It can analytically continued to whole complex plane, expect for a simple pole  at $s=1$. The value of zeta function at non-negative integers are related to the Bernoulli numbers by the following formula
\begin{align} \label{Zeta(-k)}
    \zeta(-k)=-\frac{B_{k+1}}{k+1}.
\end{align}
A family of constants appearing in the study of zeta functions and multiple gamma functions have been investigated in \cite{AdamchikNegativePolygamma,Bendersky,ChoiSomeConstants,CertainSeriesZeta,CertainSeriesGamma}. We now present the definition of so called generalized Glaisher constants or Bendersky-Adamchik constants. The family of generalized Glaisher constants $D_k$ is defined by
\begin{align*}
    \log D_k=\lim_{n \to \infty}\left( \sum_{m=1}^{n}m^k\log m -p(n,k)\right),
\end{align*}
where $p(n,k)$ is given as
\begin{align*}
    p(n,k)=\frac{n^k\log n}{2}+\frac{n^{k+1}}{k+1}\left( \log n -\frac{1}{k+1}\right)+k!\sum_{m=1}^{k}\frac{n^{k-m}B_{m+1}}{(m+1)!(k-m)!}\left( \log n+(1-\delta_{km})\sum_{l=1}^{m}\frac{1}{k-l+1} \right)
\end{align*}
and $\delta_{km}$ is the Kronecker delta symbol. Here and throughout, $k\in\mathbb{N}$ is a natural number. It can be easily seen that $D_0=(2\pi)^{1/2}$ and $D_1=A$ is the Glaisher-Kinkelin constant with the alternative definitions
\begin{align*}
    \log A&=\frac{1}{12}-\zeta'(-1)=\lim_{n \to \infty} \left( \sum_{k=1}^{n}k\log k-\left( \frac{n^2+n}{2} +\frac{1}{12}\right)\log n+\frac{n^2}{4} \right)\\
    &=\frac{1}{4}+\int_{0}^{\infty}\frac{e^{-t}}{t^2}\left( \frac{1}{e^t-1}-\frac{1}{t}+\frac{1}{2}-\frac{t}{12} \right)dt
\end{align*}
and the numerical value $A=1.2824271291\ldots$ \cite{FinchConstants}. The constants $D_2=B$ and $D_3=C$ have been studied in \cite{ChoiSomeConstants,Bendersky,CertainSeriesZeta,CertainSeriesGamma}. Similarly, expressions for constants $B$ and $C$ are obtained as 
\begin{align*}
    \log B=-\zeta'(-2)=\frac{\zeta(3)}{4\pi^2}=\lim_{n \to \infty} \left( \sum_{k=1}^{n}k^2\log k-\left(\frac{n^3}{3}+ \frac{n^2}{2} +\frac{n}{6}\right)\log n+\frac{n^3}{9}-\frac{n}{12} \right)
\end{align*}
and
\begin{align*}
    \log C=-\frac{11}{720}-\zeta'(-3)=\lim_{n \to \infty} \left( \sum_{k=1}^{n}k^3\log k-\left(\frac{n^4}{4}+\frac{n^3}{2}+ \frac{n^2}{4} -\frac{1}{120}\right)\log n+\frac{n^4}{16}-\frac{n^2}{12} \right)
\end{align*}
with their numerical values being $B=1.03091675\ldots$ and $C=0.97955746\ldots$. In fact, Adamchik \cite{AdamchikNegativePolygamma} has proved that
\begin{align} \label{D_kTheorem}
    \log D_k=\frac{B_{k+1}H_k}{k+1}-\zeta'(-k)
\end{align}
where $H_k$ are the harmonic numbers defined as
\begin{align*}
    H_k=\sum_{j=1}^k \frac{1}{j}.
\end{align*}
Furthermore, the relationship between $B$ and $\zeta(3)$ can be generalized as
\begin{align*}
    \log D_{2k}=\frac{(-1)^{k+1}(2k!)}{2(2\pi)^{2k}}\zeta(2k+1)=\frac{B_{2k}\zeta(2k+1)}{4\zeta(2k)}.
\end{align*}
Utilizing Euler-Maclaurin summation formula, another representation for generalized Glaisher constants is given in \cite{ChoiSomeConstants} as
\begin{align*}
    \log D_{k-1}=\frac{(-1)^{k+1}}{k}\sum_{m=1}^{k}\binom{k}{m}\frac{B_{k-m}}{m}+\frac{(-1)^k B_{k+1}}{k(k+1)}+\frac{(-1)^{k+1}}{k(k+1)}\int_{1}^{\infty}\frac{B_{k+1}(\left\{ x \right\})}{x^2}dx,
\end{align*}
where $\left\{ x \right\}$ denotes the fractional part of the real number $x\in \mathbb{R}$.
\newline

Now we give the definition of the Riemann-Liouville fractional integrals. Let $[a,b]$ be a finite interval on the real line, $f:[a,b]\to \mathbb{R}$ a piecewise continuous function and integrable on any subinterval of $[a,b]$. Then for $\mathfrak{Re}(\alpha)>0$, the left and right sided Riemann-Liouville fractional integrals of orders $\alpha$, $I_{a^+}^\alpha f$ and $I_{b^-}^\alpha f$ \cite{TheoryOfFractional} are respectively defined as
\begin{align*}
    I_{a^+}^\alpha f(x)=\frac{1}{\Gamma(\alpha)}\int_{a}^{x}(x-t)^{\alpha-1}f(t)dt, \qquad (b\ge x>a)
\end{align*}
and
\begin{align*}
    I_{b^-}^\alpha f(x)=\frac{1}{\Gamma(\alpha)}\int_{x}^{b}(t-x)^{\alpha-1}f(t)dt, \qquad (b>x\ge a).
\end{align*}
Taking $\alpha=1$, Riemann-Liouville fractional integrals trivially reduce to standard integration. In the special case $\alpha=n$ is a positive integer, Riemann-Liouville fractional integrals coincide with the $n$th repeated integrals
\begin{align*}
    I_{a^+}^n f(x)=\int_{a}^{x}dt_1\int_{a}^{t_1}dt_2\ldots\int_{a}^{t_{n-1}}dt_n f(t_n)=\frac{1}{(n-1)!}\int_{a}^{x}(x-t)^{n-1}f(t)dt
\end{align*}
and
\begin{align*}
    I_{b^-}^n f(x)=\int_{x}^{b}dt_1\int_{t_1}^{b}dt_2\ldots\int_{t_{n-1}}^{b}dt_nf(t_n)=\frac{1}{(n-1)!}\int_{x}^{b}(t-x)^{n-1}f(t)dt.
\end{align*}
Raabe \cite{Raabe} gave the following integral formula for the $\log \Gamma$ function
\begin{align*} 
    \int_{0}^{1}\log \Gamma(x+t)dx=\frac{\log 2\pi}{2}+t\log t-t, \quad (t>0).
\end{align*}
Upon taking $t\to0^+$, we obtain the special case
\begin{align} \label{NormalRaabe}
    \int_{0}^{1} \log \Gamma(x)dx=\frac{\log 2\pi}{2}.
\end{align}
In this paper, we obtain Raabe-type formulas for Riemann-Liouville fractional integrals of $\log \Gamma$ function. More explicitly, we consider $I_{0^+}^\alpha \left[ x^{\beta-1}\log\Gamma(x) \right](1),I_{1^-}^\alpha \left[ x^{\beta-1}\log\Gamma(x) \right](0),I_{0^+}^\alpha \left[ (1-x)^{\beta-1}\log\Gamma(x) \right](1)$ and $I_{1^-}^\alpha \left[ (1-x)^{\beta-1}\log\Gamma(x) \right](0)$ for suitable numbers $\alpha,\beta$. Our results also apply to fractional integrals involving the functions of form $x^{\beta-1}(1-x)^{\gamma-1}\log\Gamma(x)$.
\section{Main Results}
We need the following three fundamental lemmas.
\begin{lemma} \label{DigammaLemma}
    \cite{DigammaIntegralLemma} For $n\in \mathbb{N}$, the following equality holds.
    \begin{align*}
        \int_{0}^{1}x^n\psi(x)dx=\sum_{k=0}^{n-1}(-1)^k\binom{n}{k}\left( H_k\zeta(-k)+\zeta'(-k) \right)
    \end{align*}
\end{lemma}
Furthermore, substituting equations \eqref{Zeta(-k)} and \eqref{D_kTheorem} in Lemma \ref{DigammaLemma}, we obtain
\begin{align} \label{DigammaIntegral}
    \int_{0}^{1}x^n\psi(x)dx=\sum_{k=0}^{n-1}(-1)^{k+1}\binom{n}{k}\log D_k.
\end{align}
This is the form of Lemma \ref{DigammaLemma} we will use.
\begin{lemma} \label{BinomLemma}
    For $\mathfrak{Re}(\alpha)>0$ and $n\in \mathbb{N}$, the following inequality holds.
    \begin{align*}
        \left| \binom{\alpha-1}{n} \right|\le2^{\mathfrak{Re}(\alpha)-1}
    \end{align*}
\end{lemma}
\begin{proof}
    By Cauchy's residue theorem, we have the straightforward result:
    \begin{align*}
        \binom{\alpha-1}{n}=\frac{1}{2\pi i}\oint_{\left| z \right|=1}\frac{(1+z)^{\alpha-1}}{z^{n+1}}dz
    \end{align*}
    where the orientation of the circle is positive. Taking absolute values and utilizing triangle inequality, we obtain
    \begin{align*}
        \left| \binom{\alpha-1}{n} \right|=\left| \frac{1}{2\pi i}\oint_{\left| z \right|=1}\frac{(1+z)^{\alpha-1}}{z^{n+1}}dz \right|\le\frac{1}{2\pi}\oint_{\left| z \right|=1}\frac{\left| 1+z \right|^{\mathfrak{Re}(\alpha)-1}}{\left| z \right|^{n+1}}|dz|\le2^{\mathfrak{Re}(\alpha)-1}.
    \end{align*}
    Thus, the proof is complete.
\end{proof}

\begin{lemma} \label{BetaAbsConvLemma}
    For $\mathfrak{Re}(\alpha)>0$ and $\mathfrak{Re}(\beta)\ge1$, the integral
    \begin{align*}
        \int_{0}^{1}\sum_{n=0}^{\infty}(-1)^n \binom{\alpha-1}{n}t^{n+\beta-1}\log\Gamma(t)dt
    \end{align*}
    is absolutely convergent.
\end{lemma}
\begin{proof}
    For $t\in (0,1)$, it is obvious that $t^{\beta-1}\le1$. Therefore, using Lemma \ref{BinomLemma}, we obtain
    \begin{align*}
        \begin{split}
             \int_{0}^{1}\sum_{n=0}^{\infty}\left| (-1)^n \binom{\alpha-1}{n}t^{n+\beta-1}\log\Gamma(t) \right|dt&=\int_{0}^{1}\sum_{n=0}^{\infty}\left|\binom{\alpha-1}{n}\right|t^{n+\beta-1}\log\Gamma(t) dt\\
            &\le2^{\mathfrak{Re}(\alpha)-1}\int_{0}^{1}\sum_{n=0}^{\infty}t^n\log\Gamma(t) dt\\
            &=2^{\mathfrak{Re}(\alpha)-1}\int_{0}^{1}\frac{\log\Gamma(t)}{1-t}dt.
        \end{split}
    \end{align*}
    By the well-known inequality $\psi(t)\ge \log t-\frac{1}{t}$ for $t>0$, we get
    \begin{align*}
        \begin{split}
            \int_{0}^{1}\frac{\log\Gamma(t)}{1-t}dt&=-\log\Gamma(t)\log(1-t)\Big\rvert_0^1+\int_{0}^{1}\psi(t)\log(1-t)dt\\
            &=\int_{0}^{1}\psi(t)\log(1-t)dt\\
            &\le\int_{0}^{1}\log(1-t)\left( \log t-\frac{1}{t} \right)dt\\
            &=2\\
            &<\infty.
        \end{split} 
    \end{align*}
    Hence, the integral
    \begin{align*}
                \int_{0}^{1}\sum_{n=0}^{\infty}(-1)^n \binom{\alpha-1}{n}t^{n+\beta-1}\log\Gamma(t)dt
    \end{align*}
    is absolutely convergent. Thus, the proof is complete.
\end{proof}

Now we present the Raabe's formulas for $\log\Gamma$ and related functions via Riemann-Liouville fractional integrals.
\begin{theorem} \label{LeftSidedRaabe}
    For every $\mathfrak{Re}(\alpha)>0$, the following equality holds.
    \begin{align*}
        I_{0^+}^\alpha\log\Gamma(1)=\frac{1}{\Gamma(\alpha+1)}\sum_{n=1}^{\infty}(-1)^{n+1}\binom{\alpha}{n}\sum_{k=0}^{n-1}(-1)^k\binom{n}{k}\log D_k.
    \end{align*}
\end{theorem}
\begin{proof}
    By the definition of Riemann-Liouville fractional integrals, it follows that
    \begin{align*}
        I_{0^+}^\alpha\log\Gamma(1)=\frac{1}{\Gamma(\alpha)}\int_{0}^{1}(1-t)^{\alpha-1}\log\Gamma(t)dt.
    \end{align*}
    Expanding $(1-t)^{\alpha-1}$ in the integrand into binomial series yields
    \begin{align*}
        I_{0^+}^\alpha\log\Gamma(1)=\frac{1}{\Gamma(\alpha)}\int_{0}^{1}\sum_{n=0}^{\infty}(-1)^n\binom{\alpha-1}{n}t^n\log\Gamma(t)dt.
    \end{align*}
    Interchanging orders of the integral and the summation, which is permitted by Fubini-Tonelli theorems (see e.g., \cite{SteinRealAnalysis}) and Lemma \ref{BetaAbsConvLemma}, we obtain
    \begin{align} \label{AfterFubini}
        I_{0^+}^\alpha\log\Gamma(1)=\frac{1}{\Gamma(\alpha)}\sum_{n=0}^{\infty}(-1)^n\binom{\alpha-1}{n}\int_{0}^{1}t^n\log\Gamma(t)dt.
    \end{align}
    Taking $n+1$ where $n\ge0$ instead of $n$ in equation \eqref{DigammaIntegral}, we have
    \begin{align*}
        \int_{0}^{1}t^{n+1}\psi(t)dt=\sum_{k=0}^{n}(-1)^{k+1}\binom{n+1}{k}\log D_k.
    \end{align*}
    Integrating by parts, we result
    \begin{align*}
        \int_{0}^{1}t^{n+1}\psi(t)dt=t^{n+1}\log\Gamma(t)\bigg\rvert_0^1-(n+1)\int_{0}^{1}t^n\log\Gamma(t)dt=-(n+1)\int_{0}^{1}t^n\log\Gamma(t)dt
    \end{align*}
    and therefore
    \begin{align} \label{t^nlogGamma(t)Integral}
        \int_{0}^{1}t^n\log\Gamma(t)dt=\frac{1}{n+1}\sum_{k=0}^{n}(-1)^{k}\binom{n+1}{k}\log D_k.
    \end{align}
    Substituting equation \eqref{t^nlogGamma(t)Integral} into \eqref{AfterFubini} yields
    \begin{align*}
        I_{0^+}^\alpha\log\Gamma(1)=\frac{1}{\Gamma(\alpha)}\sum_{n=0}^{\infty}(-1)^n\binom{\alpha-1}{n}
        \frac{1}{n+1}\sum_{k=0}^{n}(-1)^{k}\binom{n+1}{k}\log D_k
    \end{align*}
    Now using the binomial identity $\binom{\alpha}{n+1}=\frac{\alpha}{n+1}\binom{\alpha-1}{n}$ and shifting indices in the summation, we obtain
    \begin{align*}
        I_{0^+}^\alpha\log\Gamma(1)&=\frac{1}{\Gamma(\alpha+1)}\sum_{n=0}^{\infty}(-1)^n\binom{\alpha}{n+1}\sum_{k=0}^{n}(-1)^{k}\binom{n+1}{k}\log D_k\\
        &=\frac{1}{\Gamma(\alpha+1)}\sum_{n=1}^{\infty}(-1)^{n+1}\binom{\alpha}{n}\sum_{k=0}^{n-1}(-1)^{k}\binom{n}{k}\log D_k.
    \end{align*}
    Thus, the proof is complete.
\end{proof}
\begin{corollary}
    The following formula holds: for $m\in \mathbb{N}$,
    \begin{align*}
        I_{0^+}^m\log\Gamma(1)&=\int_{0}^{1}dt_1\int_{0}^{t_1}dt_2\ldots\int_{0}^{t_{m-1}}dt_m\log\Gamma(t_m)\\
        &=\frac{1}{m!}\sum_{k=0}^{m-1}\binom{m}{k}\log D_k=\sum_{k=0}^{m-1}\frac{\log D_k}{k!(m-k)!}.
    \end{align*}
\end{corollary}
\begin{proof}
    Taking $\alpha=m\in \mathbb{N}$, interchanging the order of summations and using the identity
    \begin{align*}
                \sum_{n=k+1}^{m}(-1)^{n+k+1}\binom{m}{n}\binom{n}{k}=\frac{m}{m-k}\binom{m-1}{k},
    \end{align*}
    we obtain the desired result.
\end{proof}

\begin{remark}
    Taking $\alpha=1$ in Theorem \ref{LeftSidedRaabe}, we obtain the equation \eqref{NormalRaabe}.
\end{remark}
\begin{remark}
    Instead of Lemma \ref{DigammaLemma} and equation \eqref{t^nlogGamma(t)Integral}, we may use the result of \cite{DigammaIntegralLemma} that states for $n\in \mathbb{N}$,
    \begin{align*}
        \int_{0}^{1}t^n\log\Gamma(t)dt&=\frac{1}{n+1}\sum_{k=1}^{\left\lfloor \frac{n+1}{2} \right\rfloor}(-1)^k\binom{n+1}{2k-1}\frac{(2k)!}{2(2\pi)^{2k}}\left[ \left( \log2\pi+\gamma \right)\zeta(2k)-\zeta'(2k) \right]\\
        &+\frac{1}{n+1}\sum_{k=1}^{\left\lfloor \frac{n}{2} \right\rfloor}(-1)^{k+1}\binom{n+1}{2k}\frac{(2k)!}{2(2\pi)^{2k}}\zeta(2k+1)+\frac{\log 2\pi}{2(n+1)},
    \end{align*}
    which is essentially a more complicated version of equation \eqref{t^nlogGamma(t)Integral}, we obtain
    \begin{align} \label{OtherLeftRaabe}
        \begin{split}
        I_{0^+}^\alpha\log\Gamma(1)=\frac{1}{\Gamma(\alpha+1)}\sum_{n=1}^{\infty}(-1)^{n+1}\binom{\alpha}{n}&\left\{ \sum_{k=1}^{\left\lfloor \frac{n}{2} \right\rfloor}(-1)^k\binom{n}{2k-1}\frac{(2k)!}{2(2\pi)^{2k}}\left[\left( \log2\pi+\gamma \right)\zeta(2k)-\zeta'(2k) \right]\right.\\
        &\ \ \left.+\sum_{k=1}^{\left\lfloor \frac{n-1}{2} \right\rfloor}(-1)^{k+1}\binom{n}{2k}\frac{(2k)!}{2(2\pi)^{2k}}\zeta(2k+1)+\frac{\log 2\pi}{2} \right\}.
        \end{split}
    \end{align}
\end{remark}

\begin{theorem} \label{RightSidedRaabe}
    For every $\mathfrak{Re}(\alpha)>0$, the following equality holds.
    \begin{align*}
        I_{1^-}^\alpha\log\Gamma(0)&=\frac{1}{\Gamma(\alpha+1)}\left( \log\pi-\int_{0}^{1}\alpha t^{\alpha-1} \log\sin(\pi t)dt\right)-I_{0^+}^\alpha\log\Gamma(1)\\
        &=\frac{1}{\Gamma(\alpha+1)}\left(\sum_{n=1}^{\infty}(-1)^{n}\binom{\alpha}{n}\sum_{k=0}^{n-1}(-1)^k\binom{n}{k}\log D_k +\log\pi-\int_{0}^{1}\alpha t^{\alpha-1} \log\sin(\pi t)dt \right).  
    \end{align*}
\end{theorem}

\begin{proof}
    Starting with the definition of Riemann-Liouville fractional integrals and changing the variables to $t\to 1-t$, we have
    \begin{align} \label{LeftRL}
        I_{1^-}^\alpha\log\Gamma(0)=\frac{1}{\Gamma(\alpha)}\int_{0}^{1}t^{\alpha-1}\log\Gamma(t)dt
    \end{align}
    and
    \begin{align} \label{RightRL}
        I_{0^+}^\alpha\log\Gamma(1)=\frac{1}{\Gamma(\alpha)}\int_{0}^{1}(1-t)^{\alpha-1}\log\Gamma(t)dt=\frac{1}{\Gamma(\alpha)}\int_{0}^{1}t^{\alpha-1}\log\Gamma(1-t)dt.
    \end{align}
    Adding together equations \eqref{LeftRL} and \eqref{RightRL} and making use of the Euler's reflection formula, we get
    \begin{align}\label{Left+RightRL}
        \begin{split}
             I_{1^-}^\alpha\log\Gamma(0)+I_{0^+}^\alpha\log\Gamma(1)&=\frac{1}{\Gamma(\alpha)}\int_{0}^{1}t^{\alpha-1}\log\left( \Gamma(t)\Gamma(1-t) \right)dt\\
             &=\frac{1}{\Gamma(\alpha)}\int_{0}^{1}t^{\alpha-1}\log\frac{\pi}{\sin\pi t}dt\\
             &=\frac{1}{\Gamma(\alpha)}\left( \frac{\log\pi}{\alpha}-\int_{0}^{1}t^{\alpha-1}\log\pi t dt\right).
        \end{split}
    \end{align}
    Rearranging the equation \eqref{Left+RightRL} and substituting Theorem \ref{LeftSidedRaabe}, we obtain
    \begin{align*}
        I_{1^-}^\alpha\log\Gamma(0)&=\frac{1}{\Gamma(\alpha+1)}\left( \log\pi-\int_{0}^{1}\alpha t^{\alpha-1} \log\sin(\pi t)dt\right)-I_{0^+}^\alpha\log\Gamma(1)\\
        &=\frac{1}{\Gamma(\alpha+1)}\left(\sum_{n=1}^{\infty}(-1)^{n}\binom{\alpha}{n}\sum_{k=0}^{n-1}(-1)^k\binom{n}{k}\log D_k +\log\pi-\int_{0}^{1}\alpha t^{\alpha-1} \log\sin(\pi t)dt \right).  
    \end{align*}
    Thus, the proof is complete.
\end{proof}
\begin{corollary} \label{RightRepeatedRaabe}
    The following formula holds: for $m\in \mathbb{N}$,
    \begin{align*}
        I_{1^-}^m\log\Gamma(0)&=\int_{0}^{1}dt_1\int_{t_1}^{1}dt_2\ldots\int_{t_{m-1}}^{1}dt_m\log\Gamma(t_m)\\
        &=\frac{1}{m!}\left(\sum_{k=0}^{m-1}\binom{m}{k}\log D_k+\log\pi-\int_{0}^{1}m t^{m-1} \log\sin(\pi t)dt \right).
    \end{align*}
\end{corollary}
\begin{remark}
    Taking $\alpha=1$ in Theorem \ref{RightSidedRaabe}, we obtain the equation \eqref{NormalRaabe}.
\end{remark}
\begin{remark}
    Similarly, using the result \eqref{OtherLeftRaabe}, we have
    \begin{align*}
         I_{1^-}^\alpha\log\Gamma(0)=\frac{1}{\Gamma(\alpha+1)}\Bigg(& \log\pi-\int_{0}^{1}\alpha t^{\alpha-1} \log\sin(\pi t)dt\\ 
         &+\sum_{n=1}^{\infty}(-1)^{n}\binom{\alpha}{n}\left\{ \sum_{k=1}^{\left\lfloor \frac{n}{2} \right\rfloor}(-1)^k\binom{n}{2k-1}\frac{(2k)!}{2(2\pi)^{2k}}\left[\left( \log2\pi+\gamma \right)\zeta(2k)-\zeta'(2k) \right]\right. \\
        &\left. \left. +\sum_{k=1}^{\left\lfloor \frac{n-1}{2} \right\rfloor}(-1)^{k+1}\binom{n}{2k}\frac{(2k)!}{2(2\pi)^{2k}}\zeta(2k+1)+\frac{\log 2\pi}{2} \right\} \right).
    \end{align*}
\end{remark}
\begin{remark}
     By using the result of \cite{DigammaIntegralLemma} which states for $m\in \mathbb{N}$,
     \begin{align*}
         \int_{0}^{1}t^{m-1}\log\sin(\pi t)dt=-\frac{\log 2}{m}+(m-1)!\sum_{k=1}^{\left\lfloor \frac{m-1}{2} \right\rfloor}\frac{(-1)^k\zeta(2k+1)}{(2\pi)^{2k}(n-2k)!},
     \end{align*}
     in Corollary \ref{RightRepeatedRaabe}, we obtain
     \begin{align*}
        I_{1^-}^m\log\Gamma(0)&=\int_{0}^{1}dt_1\int_{t_1}^{1}dt_2\ldots\int_{t_{m-1}}^{1}dt_m\log\Gamma(t_m)\\
        &=\frac{1}{m!}\left(\sum_{k=0}^{m-1}\binom{m}{k}\log D_k+\log 2\pi\right)+\sum_{k=1}^{\left\lfloor \frac{m-1}{2} \right\rfloor}\frac{(-1)^{k+1}\zeta(2k+1)}{(2\pi)^{2k}(n-2k)!}.
     \end{align*}
\end{remark}

\begin{theorem} \label{BetaLeftRaabe}
        For every $\mathfrak{Re}(\alpha)>0$ and $\mathfrak{Re}(\beta)\ge1$, the following equalities hold.
        \begin{gather*}
                    I_{0^+}^\alpha \left[ x^{\beta-1}\log\Gamma(x) \right](1)=\frac{1}{\Gamma(\alpha)}\sum_{n=0}^{\infty}(-1)^n\binom{\alpha-1}{n}\Gamma(n+\beta)I_{1^-}^{n+\beta}\log\Gamma(0),
                    \\I_{1^-}^\alpha \left[ x^{\beta-1}\log\Gamma(x) \right](0)=\frac{\Gamma(\alpha+\beta-1)}{\Gamma(\alpha)}I_{1^-}^{\alpha+\beta-1}\log\Gamma(0). 
        \end{gather*}
\end{theorem}
\begin{proof}
        By the definition of Riemann-Liouville fractional integrals, it follows that
        \begin{align*}
            I_{0^+}^\alpha \left[ x^{\beta-1}\log\Gamma(x) \right](1)=\frac{1}{\Gamma(\alpha)}\int_{0}^{1}(1-t)^{\alpha-1}t^{\beta-1}\log\Gamma(t)dt.
        \end{align*}
        Expanding $(1-t)^{\alpha-1}$ into binomial series yields
        \begin{align*}
            I_{0^+}^\alpha \left[ x^{\beta-1}\log\Gamma(x) \right](1)=\frac{1}{\Gamma(\alpha)}        \int_{0}^{1}\sum_{n=0}^{\infty}(-1)^n \binom{\alpha-1}{n}t^{n+\beta-1}\log\Gamma(t)dt.
        \end{align*}
        Now we can interchange orders of the integral and the summation, which is permitted by Fubini-Tonelli theorems and Lemma \ref{BetaAbsConvLemma}, to obtain
            \begin{align*}
                \begin{split}
                    I_{0^+}^\alpha \left[ x^{\beta-1}\log\Gamma(x) \right](1)=\frac{1}{\Gamma(\alpha)}\sum_{n=0}^{\infty}(-1)^n \binom{\alpha-1}{n}\int_{0}^{1}t^{n+\beta-1}\log\Gamma(t)dt\\
                    =\frac{1}{\Gamma(\alpha)}\sum_{n=0}^{\infty}(-1)^n \binom{\alpha-1}{n}\Gamma(n+\beta)I_{1^-}^{n+\beta}\log\Gamma(0).
                \end{split}
            \end{align*}
            The second equality is trivial so we omit the proof. Thus, the proof is complete.
\end{proof}
\begin{corollary} \label{BetaLeftCorollary}
    The following formula holds: for $m\in \mathbb{N}$,
    \begin{align*}
        \begin{split}
             I_{0^+}^m \left[ x^{\beta-1}\log\Gamma(x) \right](1)&=\int_{0}^{1}dt_1\int_{t_1}^{1}dt_2\ldots\int_{t_{m-1}}^{1}t_m^{\beta-1}dt_m\log\Gamma(t_m)\\
            &=\frac{1}{(m-1)!}\sum_{n=0}^{m-1}(-1)^n\binom{m-1}{n}\Gamma(n+\beta)I_{1^-}^{n+\beta}\log\Gamma(0).
        \end{split}
    \end{align*}
\end{corollary}
\begin{corollary}
    The following formula holds: for $1\le k\in \mathbb{N}$,
    \begin{align*}
         I_{0^+}^\alpha \left[ x^{k-1}\log\Gamma(x) \right](1)=\frac{1}{\Gamma(\alpha)}\sum_{n=0}^{\infty}(-1)^n\binom{\alpha-1}{n}(n+k-1)!\ I_{1^-}^{n+k}\log\Gamma(0).
    \end{align*}
\end{corollary}

\begin{corollary}
    The following formula holds: for $m\in \mathbb{N}$ and $1\le k\in \mathbb{N}$,
    \begin{align*}
         I_{0^+}^m \left[ x^{k-1}\log\Gamma(x) \right](1)&=\int_{0}^{1}dt_1\int_{t_1}^{1}dt_2\ldots\int_{t_{m-1}}^{1}t_m^{k-1}dt_m\log\Gamma(t_m)\\
        &=\sum_{n=0}^{m-1}(-1)^n\frac{(n+k-1)\dots(n+1)}{(m-1-n)!}I_{1^-}^{n+k}\log\Gamma(0),
    \end{align*}
    where the empty product $(n+k-1)\dots(n+1)$ when $k=1$ is assumed to be $1$ .
\end{corollary}
\begin{remark}
    Taking $\beta=1$ in Theorem \ref{BetaLeftRaabe}, we obtain Theorem \ref{LeftSidedRaabe}.
\end{remark}
\begin{remark}
    Taking $\alpha=1$ and $\beta=1$ in Theorem \ref{BetaLeftRaabe}, we obtain equation \eqref{NormalRaabe}.
\end{remark}

\begin{theorem} \label{BetaRightRaabe}
        For every $\mathfrak{Re}(\beta)>0$ and $\mathfrak{Re}(\alpha)\ge1$, the following equalities hold.
        \begin{gather*}
                I_{0^+}^\alpha \left[ (1-x)^{\beta-1}\log\Gamma(x) \right](1)=\frac{\Gamma(\alpha+\beta-1)}{\Gamma(\alpha)}I_{0^+}^{\alpha+\beta-1}\log\Gamma(1), \\
                I_{1^-}^\alpha \left[ (1-x)^{\beta-1}\log\Gamma(x) \right](0)=\frac{1}{\Gamma(\alpha)}\sum_{n=0}^{\infty}(-1)^n\binom{\beta-1}{n}\Gamma(n+\alpha)I_{1^-}^{n+\alpha}\log\Gamma(0).
        \end{gather*}
\end{theorem}
\begin{proof}
    The proof is completely analogous to of Theorem \ref{BetaLeftRaabe}. Starting with the definition of Riemann-Liouville fractional integrals, it follows that
        \begin{align*}
            I_{1^-}^\alpha \left[ (1-x)^{\beta-1}\log\Gamma(x) \right](0)=\frac{1}{\Gamma(\alpha)}\int_{0}^{1}(1-t)^{\beta-1}t^{\alpha-1}\log\Gamma(t)dt.
        \end{align*}
        Expanding $(1-t)^{\beta-1}$ into binomial series yields
        \begin{align*}
            I_{1^-}^\alpha \left[ (1-x)^{\beta-1}\log\Gamma(x) \right](0)=\frac{1}{\Gamma(\alpha)}        \int_{0}^{1}\sum_{n=0}^{\infty}(-1)^n \binom{\beta-1}{n}t^{n+\alpha-1}\log\Gamma(t)dt.
        \end{align*}
        Again, we now interchange orders of the integral and the summation, which is permitted by Fubini-Tonelli theorems and Lemma \ref{BetaAbsConvLemma} (with the places of $\alpha$ and $\beta$ switched), to obtain
            \begin{align*}
                \begin{split}
                    I_{1^-}^\alpha \left[ (1-x)^{\beta-1}\log\Gamma(x) \right](0)&=\frac{1}{\Gamma(\alpha)}\sum_{n=0}^{\infty}(-1)^n \binom{\beta-1}{n}\int_{0}^{1}t^{n+\alpha-1}\log\Gamma(t)dt\\
                    &=\frac{1}{\Gamma(\alpha)}\sum_{n=0}^{\infty}(-1)^n \binom{\beta-1}{n}\Gamma(n+\alpha)I_{1^-}^{n+\alpha}\log\Gamma(0).
                \end{split}
            \end{align*}
            The first equality is trivial so we omit the proof. Thus, the proof is complete.
\end{proof}
\begin{corollary} \label{BetaRightCorollary}
    The following formula holds: for $m\in \mathbb{N}$,
    \begin{align*}
        \begin{split}
             I_{1^-}^m\left[ (1-x)^{\beta-1}\log\Gamma(x) \right](0)&=\int_{0}^{1}dt_1\int_{t_1}^{1}dt_2\ldots\int_{t_{m-1}}^{1}(1-t_m)^{\beta-1}dt_m\log\Gamma(t_m)\\
            &=\sum_{n=0}^{\infty}(-1)^n(m+n-1)\dots m \binom{\beta-1}{n} I_{1^-}^{n+m}\log\Gamma(0)
        \end{split}
    \end{align*}
    where the empty product $(m+n-1)\dots m$ when $n=0$ is assumed to be $1$ .
\end{corollary}
\begin{corollary}
    The following formula holds: for $1\le k\in \mathbb{N}$,
    \begin{align*}
         I_{1^-}^\alpha \left[ (1-x)^{k-1}\log\Gamma(x) \right](0)=\frac{1}{\Gamma(\alpha)}\sum_{n=0}^{k-1}(-1)^n\binom{k-1}{n}\Gamma(n+\alpha)I_{1^-}^{n+\alpha}\log\Gamma(0).
    \end{align*}
\end{corollary}
\begin{corollary}
    The following formula holds: for $m\in \mathbb{N}$ and $1\le k\in \mathbb{N}$,
    \begin{align*}
        \begin{split}
            I_{1^-}^m\left[ (1-x)^{k-1}\log\Gamma(x) \right](0)&=\int_{0}^{1}dt_1\int_{t_1}^{1}dt_2\ldots\int_{t_{m-1}}^{1}(1-t_m)^{k-1}dt_m\log\Gamma(t_m)\\
            &=\sum_{n=0}^{k-1}(-1)^n(m+n-1)\dots m \binom{k-1}{n} I_{1^-}^{n+m}\log\Gamma(0).
        \end{split}
    \end{align*}
\end{corollary}
\begin{remark}
    Taking $\beta=1$ in Theorem \ref{BetaRightRaabe}, we obtain Theorem \ref{RightSidedRaabe}.
\end{remark}
\begin{remark}
    Taking $\alpha=1$ and $\beta=1$ in Theorem \ref{BetaRightRaabe}, we obtain equation \eqref{NormalRaabe}.
\end{remark}

It follows by easily manipulations of integrals that the fractional integrals $I_{0^+}^\alpha \left[ \cdot \right](1),I_{1^-}^\alpha \left[ \cdot \right](0)$ of the function $x^{\beta-1}(1-x)^{\gamma-1}\log\Gamma(x)$ can be calculated by the means of Theorem \ref{BetaLeftRaabe} and \ref{BetaRightRaabe}.

\begin{remark}
    We also note that numerical verification lead us to believe that the Lemma \ref{BetaAbsConvLemma} also holds true for $\mathfrak{Re}(\beta)>0$. We have not been able to find a rigorous proof of this statement. However, assuming the Lemma \ref{BetaAbsConvLemma} also holds true for $\mathfrak{Re}(\beta)>0$, we obtain a wider range of validity for Theorem \ref{BetaLeftRaabe} and Theorem \ref{BetaRightRaabe}.
\end{remark}

\bibliography{main}
\end{document}